\documentclass[a4paper]{article}

\usepackage[english]{babel}
\usepackage{ucs}
\usepackage{url}
\usepackage[utf8x]{inputenc}
\usepackage{amsmath, amssymb, amsthm, amsfonts}
\usepackage{graphicx}
\usepackage[pdftex]{thumbpdf}
\usepackage[colorlinks,urlcolor=blue,citecolor=black,linkcolor=black,bookmarks=false]{hyperref}
\usepackage{float}
\usepackage{afterpage}
\usepackage{setspace}

\newtheorem{theorem}{Theorem}[section]
\newtheorem{lemma}{Lemma}[section]
\newtheorem{defn}{Definition}[section]
\newtheorem{assumption}{Assumption}[section]
\newtheorem{problem}{Problem}[section]
\newtheorem{proposition}{Proposition}[section]

\DeclareMathOperator{\acos}{acos}

\DeclareMathOperator{\dist}{dist}

\renewcommand{\vec}[1]{\mathbf{#1}}

\hypersetup{pdfauthor={Jorn H. Baayen and Wubbo J. Ockels},
            pdftitle={Tracking control with adaption of kites},
            pdfkeywords={kite, energy, control}}

\title{Tracking control with adaption\\
of kites\footnote{This paper is a preprint of a paper submitted to IET Control Theory \& Applications and is subject to Institution of Engineering and Technology Copyright. If accepted, the copy of record will be available at \bf{IET Digital Library}.}}

\date{June 27, 2011}

\author{
  J. H. Baayen\thanks{Independent Consultant, Weserstra\ss e 153, 12045 Berlin, Germany, jorn.baayen@gmail.com. Formerly at Delft Institute of Applied Mathematics, Faculty of Electrical Engineering, Mathematics and Computer Science, Mekelweg 4, 2628 CD, Delft, The Netherlands.} \ and
  W. J. Ockels\thanks{ASSET Institute, Kluyverweg 1, 2629 HS, Delft, The Netherlands, w.j.ockels@tudelft.nl.}
}

\begin{document}

\maketitle

\begin{abstract}
A novel tracking paradigm for flying geometric trajectories using tethered
kites is presented. It is shown how the differential-geometric notion of
turning angle can be used as a one-dimensional representation of the kite
trajectory, and how this leads to a single-input single-output (SISO) tracking
problem.  Based on this principle a Lyapunov-based nonlinear adaptive
controller is developed that only needs control derivatives of the kite
aerodynamic model. The resulting controller is validated using simulations with
a point-mass kite model.
\end{abstract}

\section{Introduction}

In 1980, Loyd wrote a seminal paper exploring the possibility of generating
electrical power using the pulling force of tethered airfoils, i.e., kites
\cite{Loyd80}. During the oil glut of the 1980s, however, interest in wind
energy dropped \cite{AWEA09}. It was only during the turn of the 21$^\textrm{st}$ century
that interest in kite power (and traction) picked up again with the work of
Meijaard, Ockels, Schwab and Diehl \cite{Meijaard99, Ockels01, DiehlPhd}.

Currently, almost the whole wind energy market makes use of horizontal axis
wind turbines. Wind turbines harvest the energy relatively close to the ground.
But realizing that up to a height of approximately $10$ kilometres the wind velocity
increases with altitude \cite{Lansdorp08}, and that the wind power raises with
the wind speed to the third power \cite{Loyd80}, one concludes that there is a
large potential to extract wind energy at higher altitudes. A research group at
the Institute for Applied Sustainable Science, Engineering and Technology
(ASSET) at Delft University of Technology is investigating a wind energy
converter using kites as lifting surfaces, i.e., the Laddermill
\cite{Ockels01}.

\begin{figure}[h]
\begin{minipage}[t]{0.47\linewidth}

\centering

\includegraphics[trim=0mm 0mm 0mm 0mm,width=160pt]{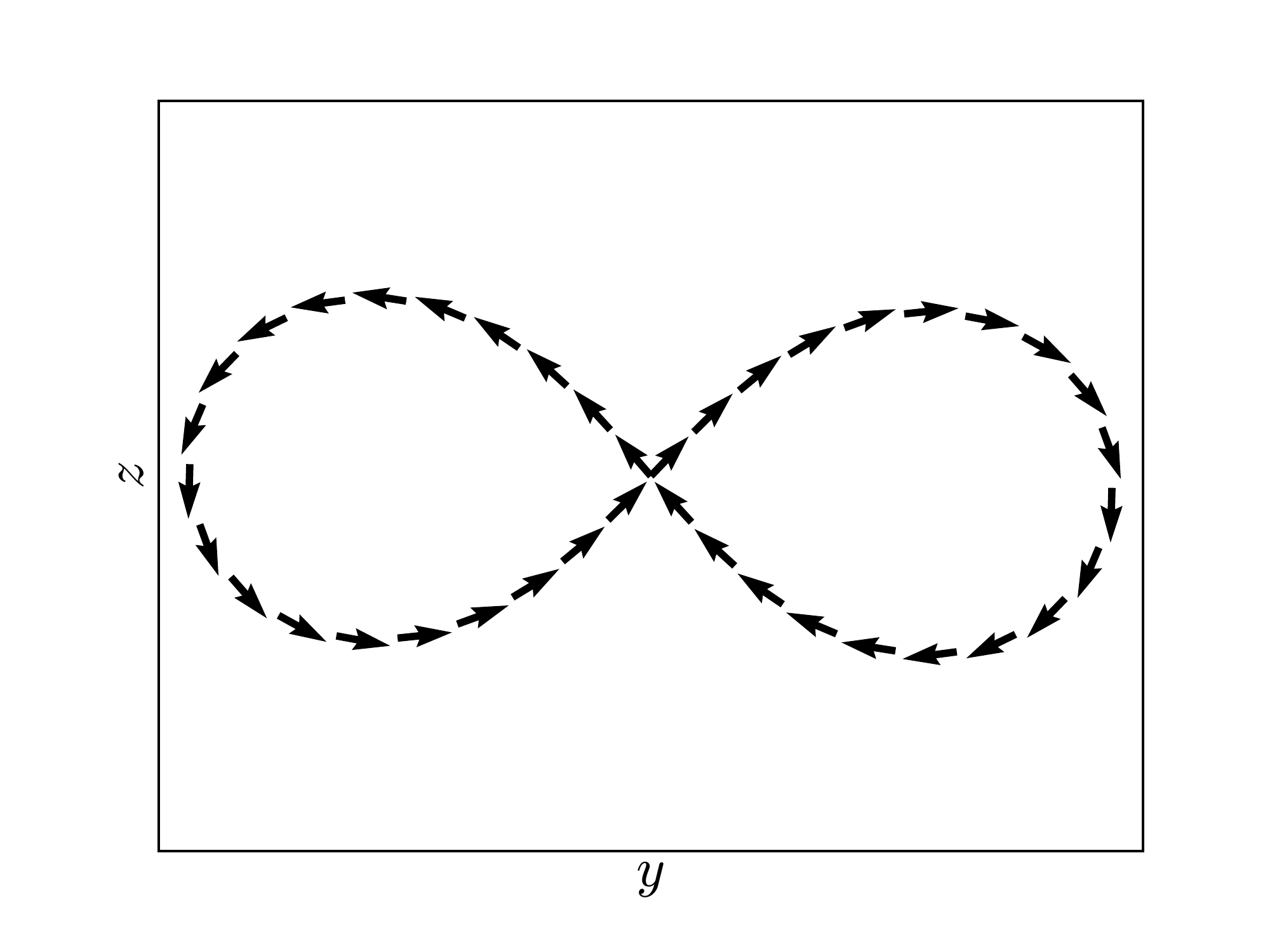}
\caption{Tangent vectors along path.}
\label{fig:figure8-tangent}

\end{minipage}
\hspace{0.3cm}
\begin{minipage}[t]{0.47\linewidth}

\centering

\includegraphics[trim=0mm 0mm 0mm 0mm,width=160pt]{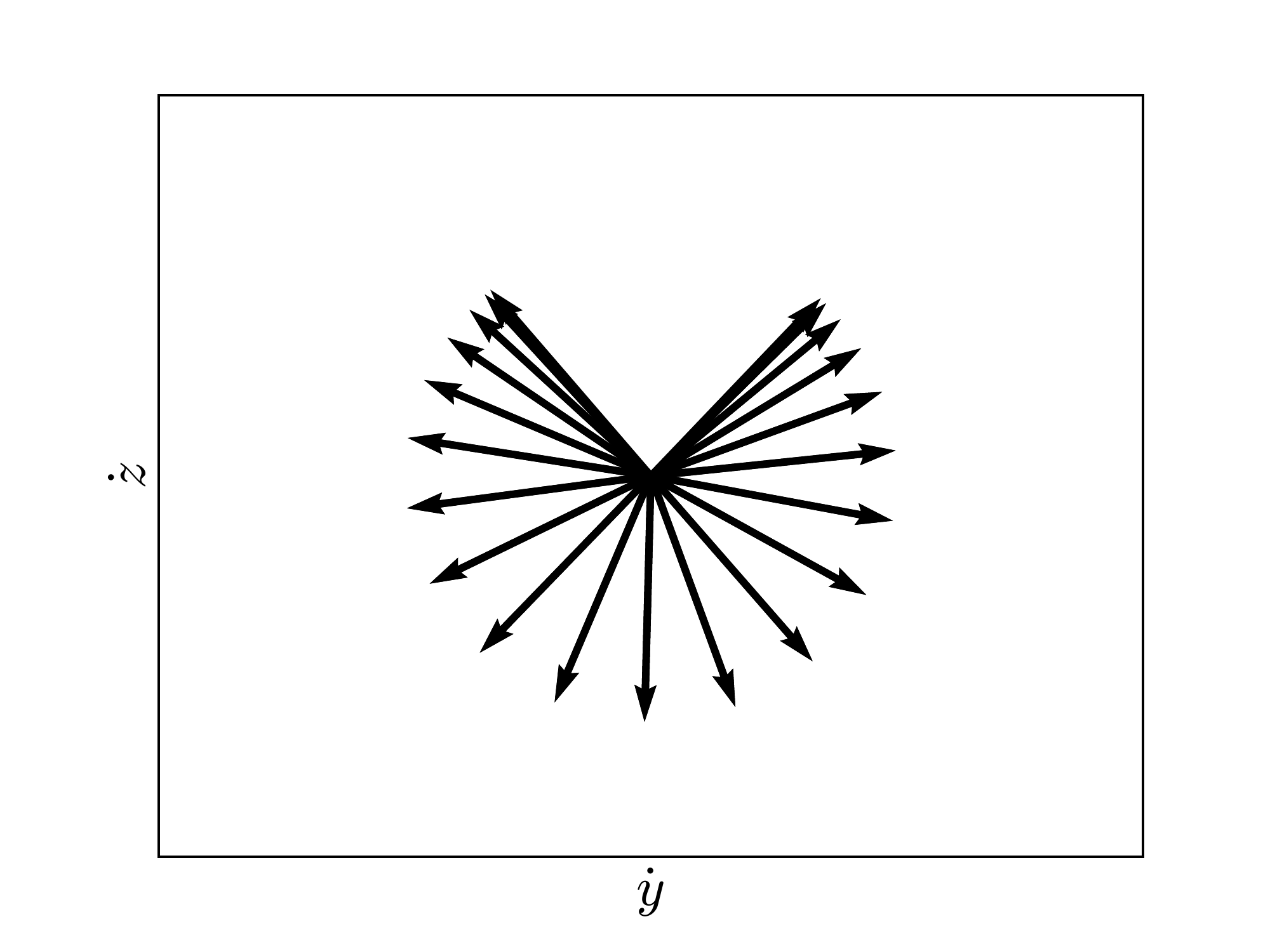}
\caption{Tangent vectors translated to the origin.}
\label{fig:figure8-tangent-angle}

\end{minipage}
\end{figure}

The design of automatic controllers for kite power systems has
been a subject of ongoing research \cite{Williams07, Williams08, Ilzhofer06,
HouskaMsc, Houska07, FagianoPhd, deGroot10}. In this paper we discuss two
contributions. First, we illustrate a trajectory tracking concept based on the
angle between the velocity vector and the horizontal. This concept is
restricted to the motion of the kite in directions perpendicular to the vector
connecting the earth tether attachment point and the kite, leaving control of
the reeling in and reeling out of the tether up to a separate controller.  To
get a first intuitive picture of this idea, let us for a moment project the
crosswind figure-eight trajectory of a kite onto the plane perpendicular to the
wind vector. As a kite traces out a figure eight on this plane, the angle the
velocity vector makes with the horizontal oscillates between certain minimum
and maximum values; see Figures \ref{fig:figure8-tangent} and
\ref{fig:figure8-tangent-angle}.

Instead of measuring this angle directly we use the more general notion of
\emph{turning angle} \cite{DifferentialGeometryBook} from elementary
differential geometry (see the Appendices).  This angle is not restricted to
$[-\pi,\pi]$ and it is continuous. It is cumulative in the sense that multiple
turns around the origin of the velocity vector will result in a multiple of
$2\pi$, but winding backwards will reduce the turning angle by the same amount
for every turn. This property makes the turning angle well suited for control
purposes. By considering the arc length as independent variable, the turning
angle -- together with initial conditions -- becomes a representation of the
image of a curve that is independent of parametrization.

Secondly, we derive an adaptive turning angle tracking controller that does not
need a full aerodynamic kite model. Instead it only uses the control
derivatives relating control input increments to steering force increments; the
rest of the required data is measured directly. We provide a stability proof
for the resulting controller and demonstrate its performance with simulations.

The remainder of this paper is organized as follows. In the next section we
discuss the physics of the kite control problem and touch on previous research.
In section 3 we formalize a geometric model of the dynamics of the kite
trajectory, and in section 4 we formulate and solve the control problem. We
illustrate the performance of the control law with simulations in section 5.
Finally, the appendices list relevant results from the field of differential
geometry.

\section{Motivation}

\begin{figure}[h]

\centering

\includegraphics[width=200pt]{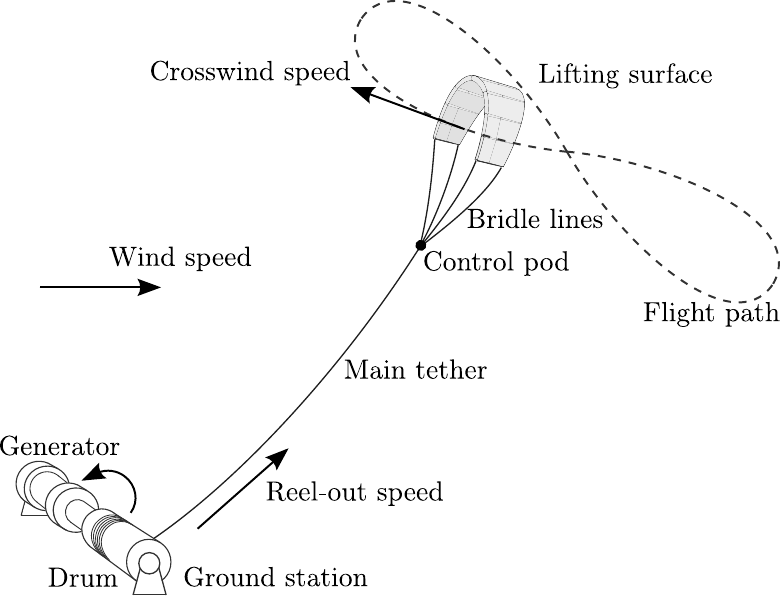}
\caption{The Laddermill during the reel-out phase \cite{deGroot10}.}
\label{fig:laddermill}

\end{figure}

The Laddermill uses kites to pull a tether from a drum driving a generator to
produce electricity. The operational mode consists of two distinct phases: a
reeling-out phase pulling the tether from the drum producing power and a
reeling-in phase which consumes energy; see Figure \ref{fig:laddermill}. It
is desired to control the kite in such a manner that a net power production
remains between the reeling-out and reeling-in phase.  By flying the kite in
directions perpendicular to the wind, the magnitude of the airflow along the
airfoil is increased \cite{Loyd80}.  This principle is comparable to the
crosswind rotation of wind turbine blades.

\begin{figure}[h]

\centering

\includegraphics[width=160pt]{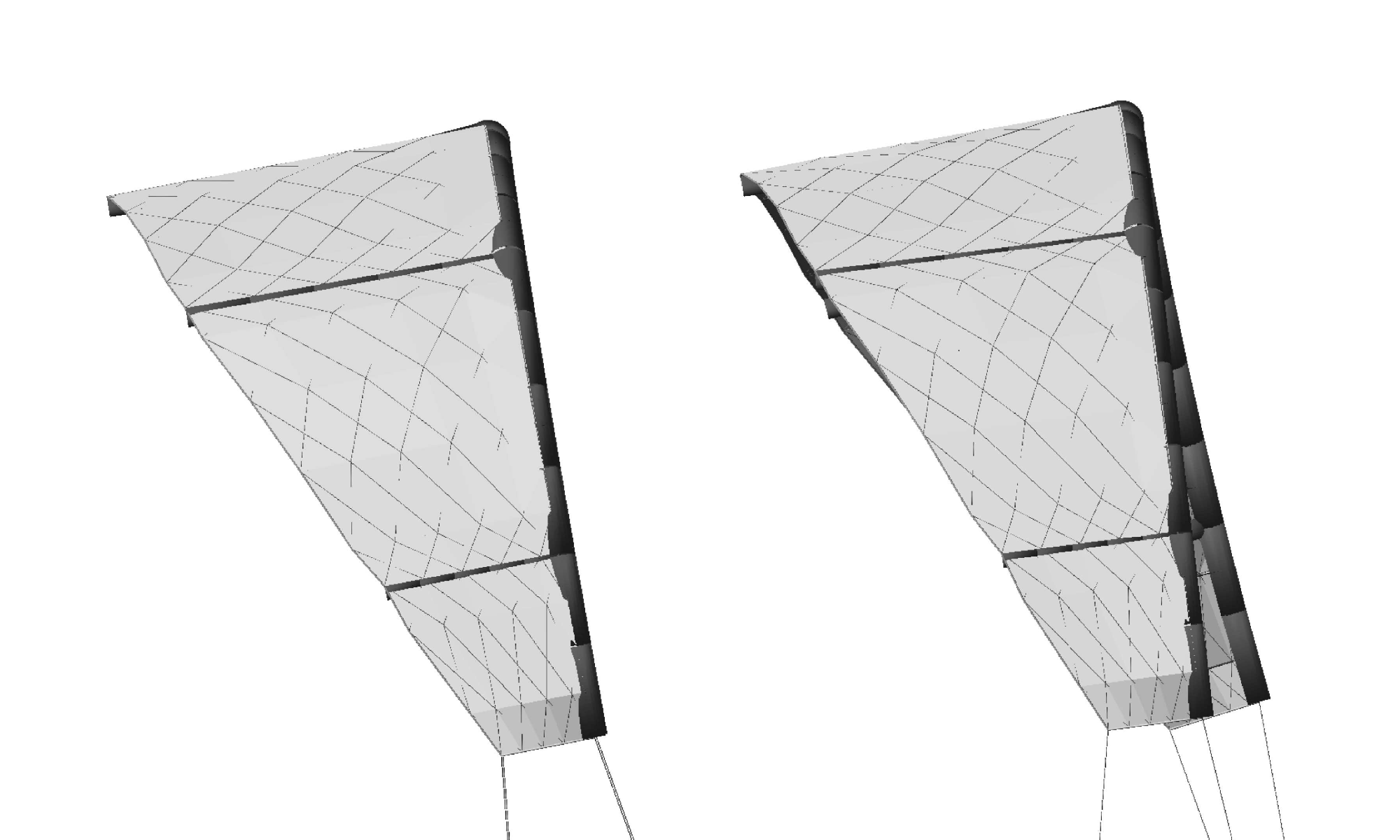}
\caption{Kite before and after applying a steering input \cite{BreukelsPhd}.}
\label{fig:kite-deformation}

\end{figure}

The control mechanism developed at the ASSET institute consists of a control
pod suspended under the kite. The control pod is connected to the drum on the
ground using a single tether and four bridle lines lead from the control pod to
the kite. The pod contains two servo motors, one controlling the length
difference of the two back lines, the \emph{steering lines}, and the other
motor controlling the overall length of the steering lines. By changing the
length difference of the steering lines the kite deforms; see Figure
\ref{fig:kite-deformation}. The aerodynamics in turn generate steering yaw
moments. A detailed explanation is available in \cite{BreukelsPhd}.  Regulation
of the total length of the steering lines provides some control over the attack
angle of the kite itself. This is known as the \emph{powering} and the
\emph{depowering} of the kite.

Previous research into automatic control of kites has focused on trajectory
tracking using linearization \cite{Williams07, Williams08} and Model Predictive
Control (MPC) \cite{Ilzhofer06, HouskaMsc, Houska07, FagianoPhd}. All of these
approaches require a fully validated kite model that can be applied in real
time. Multi-body models such as those developed in \cite{BreukelsPhd} are
precise but too computationally intensive for control design, and condensation
of kite dynamics into a rigid body model using system identification techniques
has proven to be difficult \cite{deGroot10}.

\section{Geometric modeling and system}

We are interested in controlling flight in directions perpendicular to that of
the radial vector $\vec{r}$ connecting the earth tether attachment point and
the kite. The control of the tether and of the power lines is left to a separate
controller.  In order to derive a control system that is invariant of the tether
length, we project the flight trajectory onto the upper half of the unit sphere
centered on the earth tether attachment point. This yields the projected path
$\boldsymbol\gamma:= \vec{r} / \|\vec{r}\|$.

\subsection{Geometric model}

As mentioned in the introduction, we will control the turning angle of the kite
trajectory. In order to be able to compute the turning angle on the upper half
of the unit sphere we need a coordinate system. For this, note that the
coordinate patch
\[
\vec{p}: [0, 2\pi) \times [0, \pi / 2] \to \mathbb{R}^3,\quad(v,w) \mapsto \left( \cos w \cos v, \cos w \sin v, \sin w \right).
\]
yields a continuous one-to-one (except for the pole) correspondence between
the upper half of the sphere and the set $[0, 2\pi) \times [0, \pi / 2]$. The
vector fields $\partial \vec{p} / \partial v$ and $\partial \vec{p} / \partial
w$ are orthogonal, and by normalizing them we obtain an orthonormal basis of the
tangent plane,
\begin{equation}\label{eq:tangent-basis}
\vec{e}_1 := (-\sin v,\, \cos v,\, 0), \quad \vec{e}_2 := ( -\cos v \sin w, \, -\sin v \sin w, \, \cos w),
\end{equation}
as illustrated in Figure \ref{fig:sphere}.

\begin{figure}[h]

\centering

\includegraphics[width=220pt]{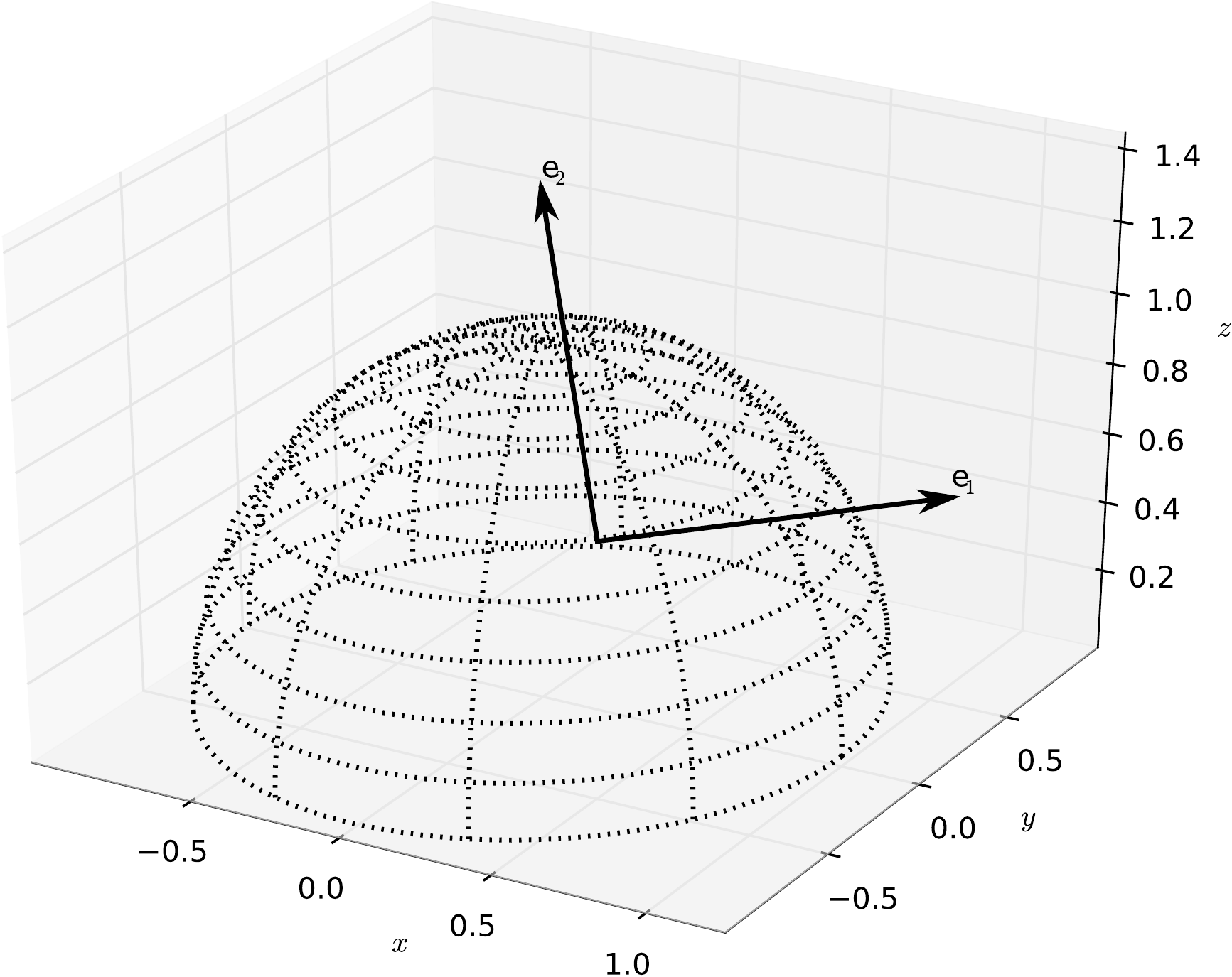}
\caption{Coordinate system on the upper half of the unit sphere.}
\label{fig:sphere}

\end{figure}

The vector field $\vec{e}_1$ has no vertical component, which intuitively makes
it a suitable reference for the definition of a turning angle. To make this
precise, we apply Liouville's theorem \cite{DifferentialGeometryBook} for which
we need the curvatures of the coordinate system. The coordinate curves $w
\mapsto \vec{p}(v,w)$ are parts of great circles, whence their geodesic
curvature $(\kappa_g)_2$ is zero. The curves $v \mapsto \vec{p}(v,w)$, however,
have geodesic curvature $(\kappa_g)_1= \sin w$.  By Liouville's theorem, the
turning angle of the projected kite trajectory with respect to the vector field
$\vec{e}_1$ is given by the integral equation
\begin{equation}\label{eq:theta}
\theta[\boldsymbol\gamma]:=\theta[\boldsymbol\gamma,\vec{e}_1] = \int_{\boldsymbol\gamma}\! \kappa_g[\boldsymbol\gamma] - (\kappa_g)_1 \cos \theta[\boldsymbol\gamma, \vec{e}_1] \,d s + \theta[\boldsymbol\gamma,\vec{e}_1]_0,
\end{equation}
where $\kappa_g[\boldsymbol\gamma]$ is the geodesic curvature of the projected
kite trajectory and $\theta[\boldsymbol\gamma,\vec{e}_1]_0$ is the initial
value of the turning angle.  The geodesic curvature of the projected kite
trajectory can be calculated as
\begin{equation}\label{eq:kappa}
\kappa_g[\boldsymbol\gamma] = \frac{\dot{x}\ddot{y} - \ddot{x}\dot{y}}{\|\dot{\boldsymbol\gamma}\|^3},
\end{equation}
with $\dot{x}$ and $\dot{y}$ the coordinates of the velocity of the projected
trajectory, expressed in the $(\vec{e}_1,\vec{e}_2)$-basis of the tangent
plane.

\subsection{Path length}\label{sec:path-length}

The projected velocity of the kite may sometimes have a component perpendicular
to the target trajectory $\boldsymbol\gamma_t$. If we zig-zag around the target
trajectory and straighten out the traversed path we obtain a longer arc length
than if we would straighten out the target trajectory.  We need to take this
difference into account, for otherwise we would move through the target trajectory
quicker than intended. For this reason we introduce a corrected path
length that disregards the perpendicular component.

Let $\vec{T}_t := \cos \theta[\boldsymbol\gamma_t] \vec{e}_1 + \sin
\theta[\boldsymbol\gamma_t] \vec{e}_2$, where $(\vec{e}_1, \vec{e}_2)$ is the
orthonormal basis of the tangent plane at $\boldsymbol\gamma_t$ given by
(\ref{eq:tangent-basis}).  This unit vector field is tangent to the
target trajectory. We transport this tangent vector to the tangent space
at $\boldsymbol\gamma$ by means of the operator $\operatorname{R}$, which
rotates $\boldsymbol\gamma$ to $\boldsymbol\gamma_t$ along their
connecting geodesic.
Taking the inner product $\dot{\boldsymbol\gamma} \cdot
(\operatorname{R}^T \vec{T}_t)$, we obtain the component of the projected velocity tangential to the
target trajectory. This motivates the following definition.

\begin{defn}
We refer to the integral
\[
s_c := \int \! \dot{\boldsymbol\gamma} \cdot (\operatorname{R}^T \vec{T}_t)\, dt
\]
as the \emph{corrected path length}.
\end{defn}

\subsection{Control system}

Our objective is to track the given target trajectory $\boldsymbol\gamma_t$
which, as hinted at previously, we will do by controlling the turning angle
$\theta[\boldsymbol\gamma]$.  We show that we can achieve this objective using a single
steering control input.  This way we construct a one-dimensional single-input,
single-output tracking problem.

In order to arrive at a set of equations that are affine in their control input
we base our discussion on steering input \emph{increments} $\tilde{u}:=u-u^*$.
Central to this linearizing approach is the following assumption.
\begin{assumption}
The acceleration $\ddot{\boldsymbol\gamma}$ is locally affine in
feasible\footnote{We call a control input increment feasible when the actuator
can implement it within a single sampling interval.} control increments
$\tilde{u}$. That is,
\[
\left|\frac{\partial^3 \boldsymbol\gamma}{\partial u \partial t^2}\right| \gg \left|\frac{\partial^{4} \boldsymbol\gamma}{\partial u^2 \partial t^2} \tilde{u}\right|.
\]
\end{assumption}

Linearizing the derivative of the turning angle in the control input $u$ we
find, using (\ref{eq:theta}), (\ref{eq:kappa}), and the definition of the
corrected path length, that
\begin{equation}
\label{eq:dthetadt}
\begin{split}
\frac{d\theta[\boldsymbol\gamma]}{dt}
                            & = \left.\frac{d\theta[\boldsymbol\gamma]}{dt}\right|_{u^*} + \left.\frac{\partial^2 \theta[\boldsymbol\gamma]}{\partial u \partial t}\right|_{u^*} \tilde{u} \\
                            & = \left.\frac{d\theta[\boldsymbol\gamma]}{ds} \frac{ds}{dt}\right|_{u^*} + \left.\frac{\partial^2 \theta[\boldsymbol\gamma]}{\partial u \partial s} \frac{ds}{dt}\right|_{u^*} \tilde{u} \\
                            & = \frac{\dot{x}\ddot{y} - \ddot{x}\dot{y}}{\|\dot{\boldsymbol\gamma}\|^2} - \boldsymbol\gamma \cdot \vec{e}_z \cos\theta[\boldsymbol\gamma] \|\dot{\boldsymbol\gamma}\| + \frac{\dot{x}\frac{\partial \ddot{y}}{\partial u} - \frac{\partial \ddot{x}}{\partial u}\dot{y}}{\|\dot{\boldsymbol\gamma}\|^2} \tilde{u},
\end{split}
\end{equation}
where $\boldsymbol\gamma \cdot \vec{e}_z$ denotes the (earth) vertical
component of the projected trajectory.  In deriving our control loop in the
next section we will make use of the direct coupling between the control
increment $\tilde{u}$ and the derivative of the turning angle.  This principle is
similar to the incremental nonlinear dynamic inversion approach for control of
aircraft as developed by Sieberling et al.  \cite{Sieberling10}.

For the derivatives $\partial \ddot{x} / \partial u$ and $\partial \ddot{y} /
\partial u$, representing the sensitivities of the accelerations to changes in the steering
input, we need a model. Let $r=\|\vec{r}\|$ and note that
\[
\ddot{\boldsymbol\gamma} = \frac{d^2}{dt^2}\left(\frac{\vec{r}}{r}\right) = \frac{\ddot{\vec{r}}}{r} - \frac{\ddot{r} \vec{r}}{r^2} - 2 \frac{\dot{r} \dot{\vec{r}}}{r^2} + 2 \frac{\dot{r}^2 \vec{r}}{r^3}.
\]
Expressed in the tangential frame the first two coordinates of the radial
vector $\vec{r}$ are zero. As a result, and due to the fact that only the
control derivatives of the accelerations are non-zero, we see that the
sensitivities, expressed in the tangential frame, are given by
\begin{equation}\label{eq:ddxydu}
\begin{split}
\frac{\partial \ddot{x}}{\partial u} & = \frac{\rho v_a^2 S_t C^t_{x,u}}{2mr}, \\
\frac{\partial \ddot{y}}{\partial u} & = \frac{\rho v_a^2 S_t C^t_{y,u}}{2mr}.
\end{split}
\end{equation}
Here $C^t_{x,u}$ and $C^t_{y,u}$ are the unitless \emph{control derivatives} expressed
in the tangential frame, $\rho$ is the density of air, $v_a$ is the airspeed,
$S_t$ is the surface area of the tips and $m$ is the combined mass of the kite
and the control pod. Let $\vec{C}_u=\left[C_{x,u}\;\; C_{y,u}\;\; C_{z,u}\right]$ denote the
control derivatives expressed in a body-fixed reference frame.

Kites deform during flight and at any point in time the deformation depends on
the form the kite had previously. Hence the deformation of the kite is
\emph{path-dependent}. In order to model the control derivatives, we assume
path-dependent deformation to have a -- for control purposes -- negligible
influence on the steering aerodynamics.

\begin{assumption}
The unitless control derivatives $\vec{C}_u$, expressed in the body-fixed
reference frame, are a function of the current steering input $u$.
\end{assumption}

Note that this assumption can be relaxed by including the dependence on additional
state variables.

We approximate the function $u \mapsto \vec{C}_u$, mapping states to control
derivatives, using B-spline networks as in \cite{Sonneveldt07}. Let $\vec{w}_i$
denote the weights and $\vec{b}$ the basis functions of the B-spline network
for the $i$th control derivative with $i \in \{x, y, z\}$, so that the $i$th
B-spline network is given by $C_{i,u} = \vec{w}_i \cdot \vec{b}$. The number of
basis functions used for the spline networks must be chosen in such a way as to
avoid problems with under- or overparametrization. This can be done using a
trial and error procedure.

By continuously measuring the position, linear velocity, linear acceleration, airspeed and
attitude of the kite, we obtain a nonlinear system of the form
\begin{equation}\label{eq:control-system}
\frac{d\theta[\boldsymbol\gamma]}{dt} = f(\theta[\boldsymbol\gamma], t) + B(t) \tilde{u}.
\end{equation}
The measurements are implicit in the time-dependence of the functions $f$ and $B$,
which are given by
\[
f=\left.\frac{d\theta[\boldsymbol\gamma]}{dt}\right|_{u^*}=\frac{\dot{x}\ddot{y} - \ddot{x}\dot{y}}{\|\dot{\boldsymbol\gamma}\|^2} - \boldsymbol\gamma \cdot \vec{e}_z \cos\theta[\boldsymbol\gamma] \|\dot{\boldsymbol\gamma}\|
\]
and
\begin{equation}\label{eq:B}
B=\left.\frac{\partial^2 \theta[\boldsymbol\gamma]}{\partial u \partial t}\right|_{u^*}=\frac{\dot{x}\frac{\partial \ddot{y}}{\partial u} - \frac{\partial \ddot{x}}{\partial u}\dot{y}}{\|\dot{\boldsymbol\gamma}\|^2}.
\end{equation}
By choosing our system in this way the need for a full aerodynamic kite model vanishes.

From (\ref{eq:dthetadt}) and (\ref{eq:ddxydu}) we see that the
function $B$ is linear in the control derivatives $C_{i,u}=\vec{w}_i \cdot
\vec{b}$, that is,
\begin{equation}\label{eq:B-lambda}
B = \sum_i \lambda_i \vec{w}_i \cdot \vec{b},
\end{equation}
for some coefficients $\lambda_i \in \mathbb{R}$. This shows that the system
(\ref{eq:control-system}) is affine in both the control input $\tilde{u}$ as well
as in the control derivative spline weights $\vec{w}_i$, a property that will
allow us to design an adaptive control loop with relative ease in the next
section.

\section{Control}

In this section we design our control loop and control derivative estimators by
means of a control-Lyapunov approach. Our approach is closely related to those
described in \cite{Sonneveldt07, AdaptiveControlBook, NonlinearSystemsBook,
NonlinearAndAdaptiveControlDesignBook, FlightControlDesignBook}.

Our objective is to track a given reference projected kite trajectory. To this
end, we construct a flight controller consisting of an inner and an outer loop. The inner loop controls
the turning angle of the kite trajectory $\theta[\boldsymbol\gamma]$ towards
a target turning angle $\theta_t$. The outer loop computes a target turning angle $\theta_t$ that steers
the kite onto the target trajectory.

\subsection{Inner loop}


Our inner loop cost function is the difference between the actual and desired turning angle:
\begin{equation}\label{TrackingError}
e:=\theta[\boldsymbol\gamma] - \theta_t.
\end{equation}
In order to track the reference trajectory we must find a control
law that renders $e \to 0$.

The steering dynamics of a kite may subtly change over time due to weather
patterns and wear.  In order to automatically adapt to these changes we aim for
an adaptive control approach that maintains and updates an internal
representation of the steering control derivatives. We will denote this
internal representation as $\hat{\vec{w}}_i$ with $i \in \{x, y, z\}$.

Summarizing, we formulate the following control problem:

\begin{problem}
\label{prob:control}
Find a steering control law and a control derivative B-spline weight estimator
update laws that such that $e \to 0$ as $t \to \infty$.
\end{problem}

\subsubsection{Control design}

In this section we propose steering control and estimator update laws to
solve Problem \ref{prob:control} using a control-Lyapunov approach.

\begin{proposition}
Consider the system (\ref{eq:control-system}).
The control law
\begin{equation}\label{eq:control-law}
\tilde{u} := -\frac{1}{\sum_i \lambda_i \hat{\vec{w}}_i \cdot \vec{b}} \left[\left.\frac{d\theta[\boldsymbol\gamma]}{dt}\right|_{u^*} - \frac{d\theta_t}{dt} + Ke\right],
\end{equation}
with gain $K \in \mathbb{R}^{+}$, and the control derivative B-spline weight
estimator update laws
\begin{equation}\label{eq:update-law}
\frac{d\hat{\vec{w}}_i}{dt}:=\Gamma e \lambda_i \tilde{u} \vec{b},
\end{equation}
with adaptation gain $\Gamma \in \mathbb{R}^{+}$, render the Lyapunov function
\begin{equation}\label{eq:V}
V := \frac{1}{2}e^2 + \frac{1}{2} \Gamma^{-1} \sum_i \|\hat{\vec{w}_i} - \vec{w}_i\|^2
\end{equation}
decreasing along the flow.
\end{proposition}

\begin{proof}
The time-derivative of the Lyapunov function $V$ is given by
\begin{equation}\label{eq:dVdt}
\frac{dV}{dt} = e \frac{de}{dt}  + \Gamma^{-1} \sum_i (\hat{\vec{w}}_i - \vec{w}_i) \cdot \frac{d\hat{\vec{w}}_i}{dt}
\end{equation}
with, using (\ref{eq:dthetadt}),
\[
\frac{de}{dt} = \left.\frac{d\theta[\boldsymbol\gamma]}{dt}\right|_{u^*} + \left.\frac{\partial^2 \theta[\boldsymbol\gamma]}{\partial u \partial t}\right|_{u^*} \tilde{u} - \frac{d\theta_t}{dt}.
\]
The estimator update laws (\ref{eq:update-law}) cancel out all terms with $\vec{w}_i$.
Indeed, replacing these into (\ref{eq:dVdt}), we obtain, using (\ref{eq:B}) and (\ref{eq:B-lambda}), that
\begin{equation}\label{eq:dVdt2}
\frac{dV}{dt} = e \left[\left.\frac{d\theta[\boldsymbol\gamma]}{dt}\right|_{u^*} - \frac{d\theta_t}{dt}\right] + e \tilde{u} \sum_i \lambda_i \hat{\vec{w}}_i \cdot \vec{b}.
\end{equation}
Finally, substituting the control law (\ref{eq:control-law}) into (\ref{eq:dVdt2}), we see that
\[
\frac{dV}{dt} = -Ke^2.
\]
The gain $K \in \mathbb{R}^+$ renders $V$ strictly decreasing whenever $e \neq
0$.
\end{proof}

In the previous analysis we neglected rate and magnitude constraints of the
control signal. Let $\tilde{u}_a$ denote the control increment
implemented by the actuator.  Since in general $\tilde{u}_a \neq \tilde{u}$, the Lyapunov
function (\ref{eq:V}) is no longer guaranteed to be decreasing and the
parameter estimation process may diverge. In order to remedy this situation we
introduce the \emph{modified tracking error} as in \cite{Farrell03},
\[
e_m := e - \zeta,
\]
with
\begin{equation}\label{eq:dzetadt}
\frac{d{\zeta}}{dt} := -K \zeta + \left.\frac{\partial^2 \theta[\boldsymbol\gamma]}{\partial u \partial t}\right|_{u^*,\hat{\vec{w}}_i} (\tilde{u}_a - \tilde{u}),
\end{equation}
providing a correction for the influence of the actuator limitations, where the second
derivative is evaluated using the estimated weights $\hat{\vec{w}}_i$.  Note
that in the absence of actuator constraints $e_m \to e$.

In order to render the modified tracking error to zero, we adjust the estimator
update laws (\ref{eq:update-law}) to be proportional to $e_m$.  The
associated modified Lyapunov function is -- by construction -- decreasing along
the flow, independently of whether or not $\tilde{u}_a=\tilde{u}$:

\begin{proposition}
\label{prop:mod-V}
Consider the system (\ref{eq:control-system}).
Any control law
and the modified control derivative B-spline weight estimator update laws
\begin{equation}\label{eq:mod-update-law}
\frac{d\hat{\vec{w}}_i}{dt} := \Gamma e_m \lambda_i \tilde{u}_a \vec{b}
\end{equation}
render the modified Lyapunov function
\begin{equation}\label{eq:mod-V}
V_m := \frac{1}{2} e_m^2 + \frac{1}{2} \Gamma^{-1} \sum_i \|\hat{\vec{w}}_i - \vec{w}_i\|^2
\end{equation}
decreasing along the flow.
\end{proposition}

\begin{proof}
The time-derivative of the modified Lyapunov function (\ref{eq:mod-V}) is given by
\begin{equation}\label{eq:dV_mdt}
\frac{dV_m}{dt} = e_m \frac{de_m}{dt} + \Gamma^{-1} \sum_i (\hat{\vec{w}}_i - \vec{w}_i) \cdot \frac{d\hat{\vec{w}}_i}{dt}
\end{equation}
with, using (\ref{eq:dzetadt}) and (\ref{eq:dthetadt}),
\[
\frac{de_m}{dt} = K (e - e_m) - \left.\frac{\partial^2 \theta[\boldsymbol\gamma]}{\partial u \partial t}\right|_{u^*, \hat{\vec{w}}_i}(\tilde{u}_a - \tilde{u})  + \left.\frac{d\theta[\boldsymbol\gamma]}{dt}\right|_{u^*} + \left.\frac{\partial^2 \theta[\boldsymbol\gamma]}{\partial u \partial t}\right|_{u^*} \tilde{u}_a -
                 \frac{d\theta_t}{dt}.
\]
The modified estimator update laws (\ref{eq:mod-update-law}) cancel out all terms with $\vec{w}_i$.
Indeed, replacing these into (\ref{eq:dV_mdt}) we obtain, using (\ref{eq:B}) and (\ref{eq:B-lambda}), that
\begin{equation}\label{eq:dV_mdt2}
\begin{split}
\frac{dV_m}{dt} & = e_m \left[K (e - e_m) - \left.\frac{\partial^2 \theta[\boldsymbol\gamma]}{\partial u \partial t}\right|_{u^*, \hat{\vec{w}}_i}(\tilde{u}_a - \tilde{u})\right] \\
                 & \qquad + e_m \left[ \left.\frac{d\theta[\boldsymbol\gamma]}{dt}\right|_{u^*}  -
                 \frac{d\theta_t}{dt} \right] + e_m \tilde{u}_a \sum_i \lambda_i \hat{\vec{w}_i} \cdot \vec{b}.
\end{split}
\end{equation}
Noting that $\tilde{u}_a=(\tilde{u}_a - \tilde{u}) + \tilde{u}$ and
substituting the canonical control law (\ref{eq:control-law})
into (\ref{eq:dV_mdt2}) we see that
\[
\frac{dV_m}{dt} = -Ke_m^2.
\]
The gain $K \in \mathbb{R}^{+}$ renders $V_m$ strictly decreasing whenever $e_m \neq 0$.
\end{proof}

\subsubsection{Control performance}

In this section we investigate the performance of the proposed control and
estimator update laws.  In the case that $\tilde{u}_a=\tilde{u}$, the function
$V$ defined in (\ref{eq:V}) is decreasing along the flow. This implies
convergence of the kite trajectory turning angle $\theta[\gamma]$ to the target
turning angle $\theta_t$.

\begin{theorem}
Consider the system (\ref{eq:control-system}), the control law
(\ref{eq:control-law}) and the estimator update laws (\ref{eq:mod-update-law}).
Assume that there is a $t_0 > 0$ such that for all $t > t_0$
\[
\|\dot{\boldsymbol\gamma}(t)\| > 0
\]
holds, so that for $t > t_0$ the kite trajectory is a regular curve.
Additionally assume that for all $t > t_0$ the geometric condition
\[
\left.\frac{d \ddot{\boldsymbol\gamma}}{du}\right|_t \not\parallel \vec{r}(t),
\]
i.e. that the given two vectors are not parallel, is satisfied, so that the control law is well-defined. Finally, assume that for some
$M \in \mathbb{R}^+$ and for all $t > 0$ the bounds
\[
\left|\left.\frac{d\theta[\boldsymbol\gamma]}{dt}\right|_t\right| \leq M,\quad \left|\left.\frac{d\theta_t}{dt}\right|_t\right| \leq M, \quad |\theta[\boldsymbol\gamma](t) - \theta_t(t)| \leq M
\]
hold. Then, for any initial condition and for any reference
trajectory the turning angle of the kite trajectory converges to the target:
\[
\lim_{t \rightarrow \infty} \theta[\boldsymbol\gamma](t) = \theta_t(t).
\]
\end{theorem}

\begin{proof}
Barbalat's Lemma \cite{NonlinearSystemsBook} applied to the Lyapunov
function $V$ shows that $e \to 0$.
\end{proof}

We proceed to show that, even when $\tilde{u}_a \neq \tilde{u}$, the modified
parameter estimation process is stable in the sense that the weight estimation
errors remain bounded:

\begin{theorem}
Consider the system (\ref{eq:control-system}) and the estimator update laws
(\ref{eq:mod-update-law}). Assume that there is a $t_0 > 0$ such that
for all $t > t_0$ the conditions
\[
\|\dot{\boldsymbol\gamma}(t)\| > 0,
\]
and
\[
\left.\frac{d \ddot{\boldsymbol\gamma}}{du}\right|_t \not\parallel \vec{r}(t)
\]
are satisfied. Then for any initial estimates, any reference trajectory
and any actuated control increments $\tilde{u}_a$, the estimator errors
\[
\|\hat{\vec{w}}_i(t) - \vec{w}_i(t)\|,\quad i \in \{x, y, z\}
\]
are bounded, and the bound decreases in time for $t > t_0$.
\end{theorem}

\begin{proof}
Boundedness follows from the definition of the Lyapunov function $V_m$:
\[
\|\hat{\vec{w}}_i(t) - \vec{w}_i(t)\| \leq \sqrt{2 \Gamma V_m(t)}.
\]
Proposition \ref{prop:mod-V} shows that the bound decreases in time.
\end{proof}

We do not, however, have a guarantee that the estimates will converge to the
true values of the weights. The modified Lyapunov function $V_m$ becomes
stationary once $e_m=0$, and at this point the estimation error will not be able
to decrease further.

\subsection{Outer loop}

The task of the outer loop controller is to compute an appropriate turning angle.
To this end we consider the outer loop nonlinear system of the form
\begin{equation}\label{eq:outer-loop-system}
\frac{d \boldsymbol\gamma}{dt} = v(\cos \theta_t \vec{e}_1 + \sin \theta_t \vec{e}_2),
\end{equation}
where the velocity $v=v(t) \geq 0$ is determined by the dynamics of the kite. For this system
we obtain the following control problem:

\begin{problem}
\label{prob:outer-loop-control}
Find a turning angle control law such that $\boldsymbol\gamma \to \boldsymbol\gamma_t$ as $t \to \infty$.
\end{problem}

\subsubsection{Control design}

In this section we propose a turning angle control law to solve
Problem \ref{prob:outer-loop-control} using a control-Lyapunov approach.

Recall that $\vec{T}_t$ is a unit vector tangent to the target trajectory; cf Section \ref{sec:path-length}. Let $\vec{T}:=\operatorname{R}^T \vec{T}_t$, where the operator $\operatorname{R}$ rotates $\gamma$ to $\gamma_t$ along their connecting geodesic, so that $\vec{T}$ is an element of the tangent plane at $\boldsymbol\gamma$.
Let $\vec{t}:=\boldsymbol\gamma \times \vec{T}$, so that the vectors $\vec{T}$
and $\vec{t}$ form a basis of the tangent plane at $\boldsymbol\gamma$. Finally, denote
\[
\boldsymbol\gamma_t^{\perp}:=\boldsymbol\gamma_t-(\boldsymbol\gamma_t\cdot\boldsymbol\gamma)\boldsymbol\gamma.
\]

\begin{proposition}\label{thm:W-lyapunov}
The turning angle $\theta_t$ determined by
\begin{equation}\label{eq:target-velocity}
\cos \theta_t \vec{e}_1 + \sin \theta_t \vec{e}_2 = \vec{T} + L (\vec{t} \cdot \boldsymbol\gamma_t^{\perp}) \vec{t},
\end{equation}
with gain $L \in \mathbb{R}^+$, renders the geodesic distance
\[
W=\dist(\boldsymbol\gamma, \boldsymbol\gamma_t)
\]
decreasing along the flow.
\end{proposition}

\begin{proof}
\label{prop:outer-loop-decreasing}
The time-derivative of the geodesic distance is given by Lemma \ref{lemma:geodist-deriv} in Appendix B:
\begin{equation}\label{eq:dWdt}
\frac{dW}{dt} = -(\dot{\boldsymbol\gamma} - \operatorname{R}^T \dot{\boldsymbol\gamma_t}) \cdot \vec{Y}_{\boldsymbol\gamma}^{\boldsymbol\gamma_t},
\end{equation}
where $\vec{Y}_{\boldsymbol\gamma}^{\boldsymbol\gamma_t}$ is the geodesic vector pointing from $\boldsymbol\gamma$ towards $\boldsymbol\gamma_t$.
Substituting (\ref{eq:target-velocity}) into (\ref{eq:dWdt}), and noting that
\[
\dot{\boldsymbol\gamma}_t=(\dot{\boldsymbol\gamma} \cdot \vec{T}) \vec{T}_t,
\]
we obtain
\[
\frac{dW}{dt}=- L v (\vec{t} \cdot \vec{Y}_{\boldsymbol\gamma}^{\boldsymbol\gamma_t}) (\vec{t} \cdot \boldsymbol\gamma_t^{\perp}) = - L v \frac{\left(\vec{t} \cdot \boldsymbol\gamma_t^{\perp}\right)^2}{\|\boldsymbol\gamma_t^{\perp}\|}.
\]
Note that $dW / dt \to 0$ as $\boldsymbol\gamma_t^{\perp} \to 0$.
\end{proof}

\subsubsection{Control performance}

Proposition \ref{prop:outer-loop-decreasing} shows that the control law (\ref{eq:target-velocity}) renders the geodesic distance $\dist(\boldsymbol\gamma, \boldsymbol\gamma_t)$ decreasing
along the flow. In this sense, the control law (\ref{eq:target-velocity}) renders the target trajectory stable, but not asymptotically so. For the latter we would also need to be able
to control the kite velocity $v$.

\section{Simulation}

In this section we report on the performance of our controller with an
extension of a point-mass kite model commonly found in the literature.

\subsection{Model}

We present an extension of the point-mass kite model discussed in
\cite{DiehlPhd, HouskaMsc, Houska07, Ilzhofer06, FagianoPhd}.  For simplicity,
our presentation uses a centripetal force to keep the kite on a sphere
with given radius.

We express the aerodynamic forces acting on the kite in an aerodynamic
reference frame. Let $\vec{v}_a=\dot{\vec{r}}-\vec{v}_w$ denote the
apparent wind vector, where $\vec{v}_w$ represents the environmental wind.  Let
$\vec{e}_x$ be the unit vector collinear with the apparent wind vector
$\vec{v}_a$.  We assume that the kite cannot roll with respect to the tether,
and so we must have $\vec{e}_y$ such that $\vec{e}_y$ is parallel to
$\vec{r}~\times~\vec{e}_x$. Choose $\vec{e}_y$ to be collinear with $\vec{r}
\times \vec{e}_x$ and let $\vec{e}_z := \vec{e}_x \times \vec{e}_y$.

In expressing the magnitude of the forces we will make use of unitless
lift coefficients. Let $C_L: \alpha \mapsto C_L(\alpha)$ be the unitless
lift curve with $\alpha:=\acos(\vec{r} \cdot \vec{e}_z / \|\vec{r}\|) + \theta_p$, where
$\theta_p$ is the power setting, and let $C_S: u \mapsto C_S(u)$ be the
unitless steering force curve.  We model the drag force as the sum of a
zero-lift drag coefficient and a parabolic approximation of induced drag
\cite{AerodynamicsBook}:
\[
C_D = C_{D_0} + k C_L^2.
\]
The constant $k$ is defined by $k^{-1} = \pi A e$, with aspect ratio $A = b^2 /
S$, Oswald factor $e$, projected surface area $S$, and wingspan $b$. It follows
that the aerodynamic forces are given by
\begin{align*}
F_x & = -\frac{1}{2} \rho \|\vec{v}_a\|^2 S C_D(\alpha), \\
F_y & = -\frac{1}{2} \rho \|\vec{v}_a\|^2 S C_S(u), \\
F_z & = -\frac{1}{2} \rho \|\vec{v}_a\|^2 S C_L(\alpha),
\end{align*}
where $\rho$ is the density of air.

In our model the steering control input $u$ causes the kite to accelerate in
the direction of the Y-axis, which in turn causes an increase in the relevant
component of the apparent wind vector $\vec{v}_a$. The aerodynamic reference
frame rotates along with the apparent wind vector, and in this way the effects
of the steering yaw moments observed in real kites are emulated.  The model is
completed by considering gravity and constraint forces to restrict the kite to
a sphere.

\begin{figure}[h]

  \centering
  \includegraphics[width=320pt]{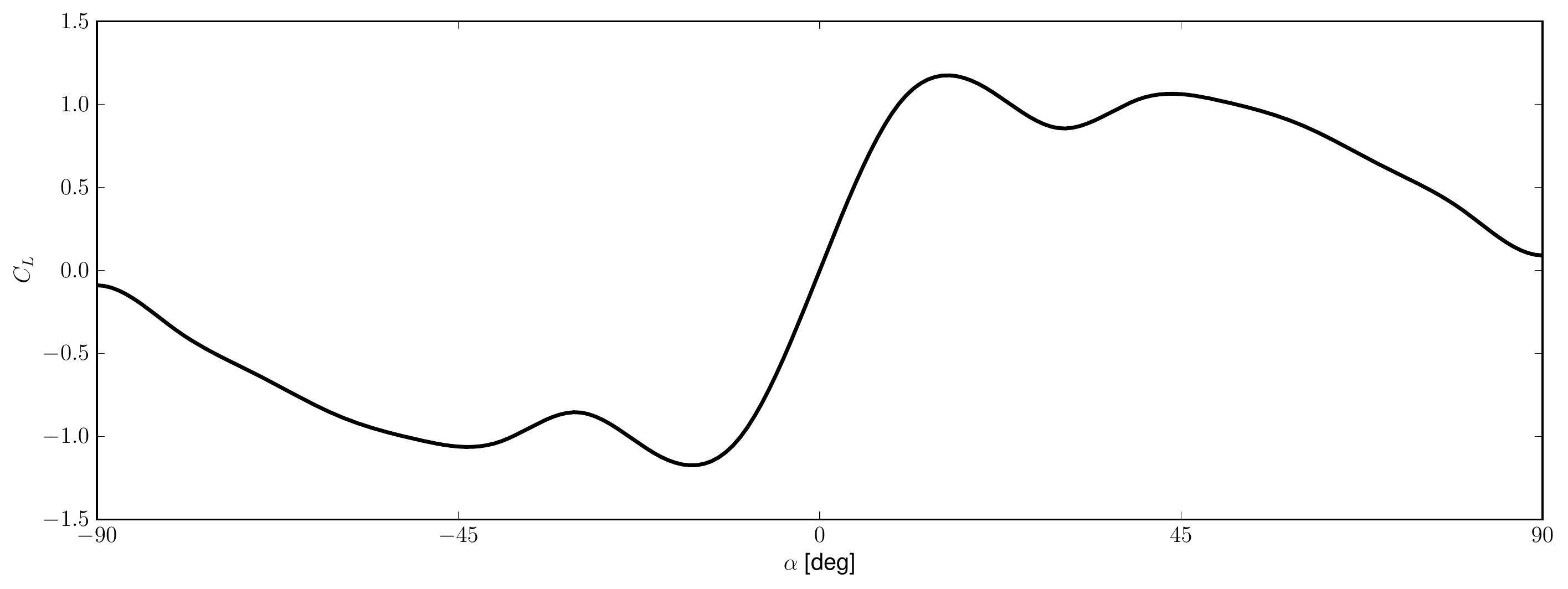}
  \caption{Lift curve modeled on data from \cite{AeroCharVertAxisWindTurb} for a NACA 0018 airfoil at Re = $2 \times 10^6$.}
  \label{fig:lift-curve}

\end{figure}

In our simulations we choose a projected surface area of $S=11$ m$^2$, a wing
span of $b=6$ m, an Oswald factor of $e=0.7$ and a mass of $m=1.5$ kg. The
zero-lift drag coefficient is given by $C_{D_0}=0.075$. The lift curve $C_L$,
modeled on data from \cite{AeroCharVertAxisWindTurb} for a NACA 0018 airfoil in a flow with Reynolds number $2 \times 10^6$ using B-splines,
is given in Figure \ref{fig:lift-curve}. Due to the current lack of more
specific data, we also set $C_S=C_L$. This is justified by the observation that one can view
the tips of the kite as vertically positioned airfoils, with a local
attack angle controlled using the steering control input $u$.

\subsection{Reference trajectory}

For the tracking objective to be independent of the velocity we treat the
corrected path length -- as opposed to time -- as the independent variable.
Our reference turning angle is
\begin{equation}\label{eq:theta_gamma_t}
\theta[\boldsymbol\gamma_t](s_{c}) := A \left[\cos\left(\frac{2\pi}{L}
s_{c}\right) - 1\right] + 0.834029,
\end{equation}
where $A \approx 2.40483$ is a root of the 0th Bessel function of the first kind
and $L=4 / 3$ is the length of the projected target trajectory.  This turning
angle generates a figure eight trajectory. It can be shown that, for the planar
notion of turning angle, the constant $A$ being a root of the appropriate
Bessel function is a necessary and sufficient condition for the turning angle
(\ref{eq:theta_gamma_t}) to result in an $L$-periodic trajectory.
The second constant is the initial condition which also corrects for the
curvature of the sphere. It was found using a root-finding procedure. The
resulting trajectory is an idealization of trajectories flown using the
prototype Laddermill kite power plant.

\subsection{Results}

The wind vector points in the direction of the earth X-axis with magnitude $6$~m/s.
The kite is tethered on a $100$ m line, and the controller gains are tuned --
using a trial and error procedure -- to $K=3.0$, $L=10.0$, and $\Gamma=100.0$.
The steering rate is limited to $0.05$ rad/s, the magnitude of the steering input
to $0.12$ rad, the power setting is fixed at $\theta_p=0.01$ rad, and every control
derivative estimator B-spline has $10$ knots.

\begin{figure}[h]

\centering

\includegraphics[width=220pt]{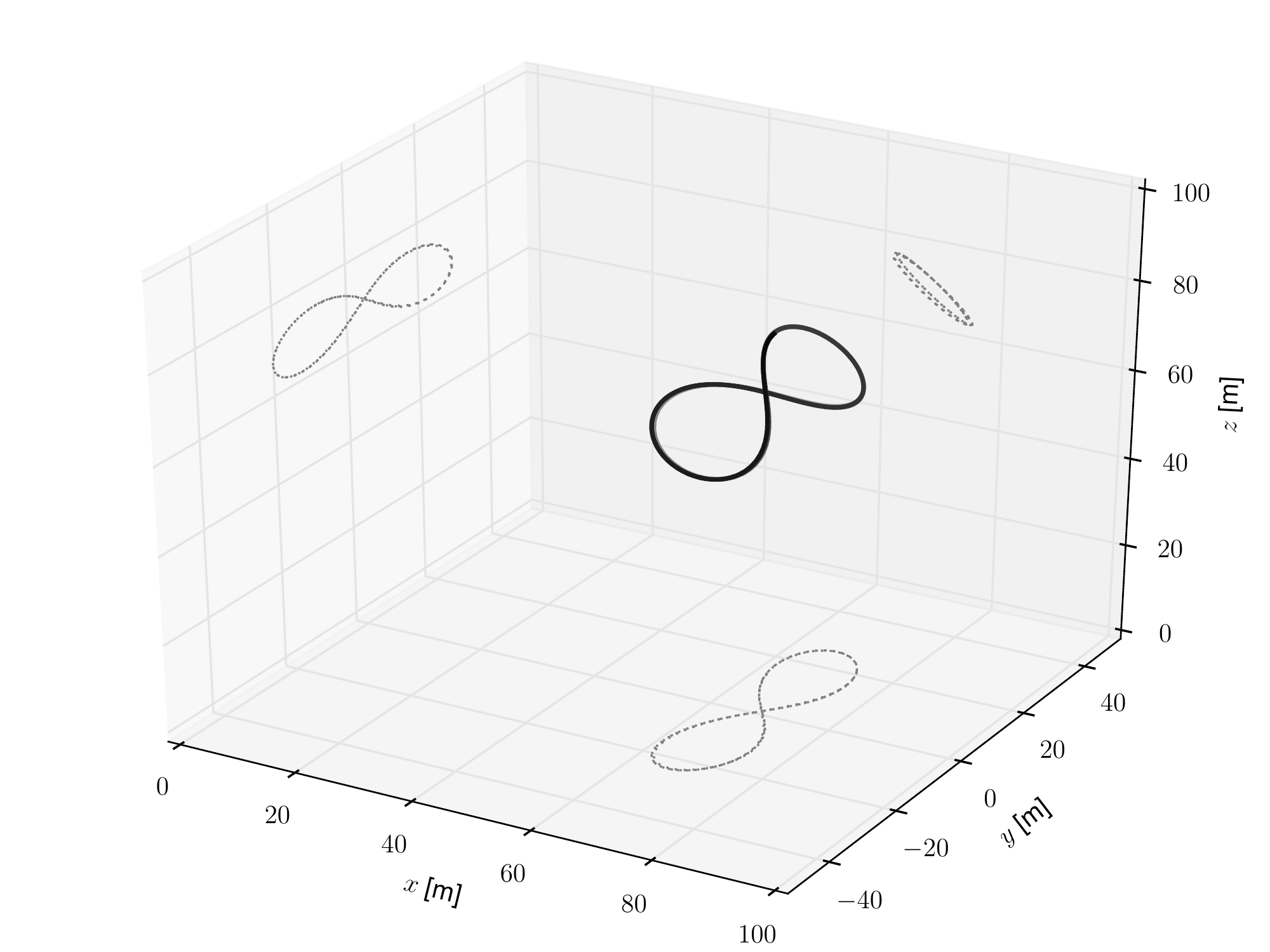}
\caption{Three overlapping figure eight trajectories.}
\label{fig:trajectory}

\end{figure}

\begin{figure}[h]

\begin{minipage}[t]{0.47\linewidth}

\centering

\includegraphics[width=160pt]{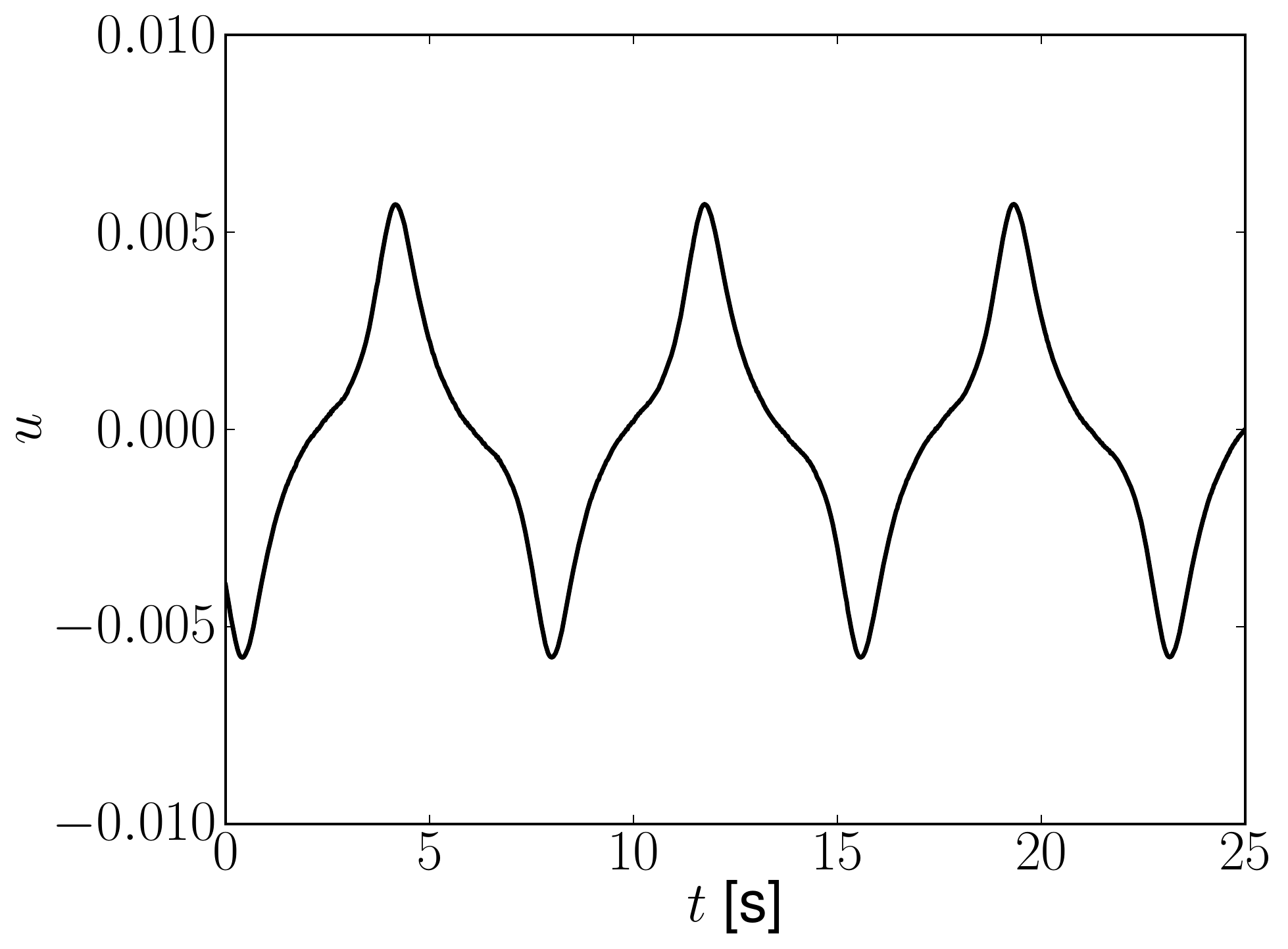}
\caption{Control input for the trajectory shown in Figure \ref{fig:trajectory}.}
\label{fig:steering}

\end{minipage}
\hspace{0.3cm}
\begin{minipage}[t]{0.47\linewidth}

\centering

\includegraphics[width=160pt]{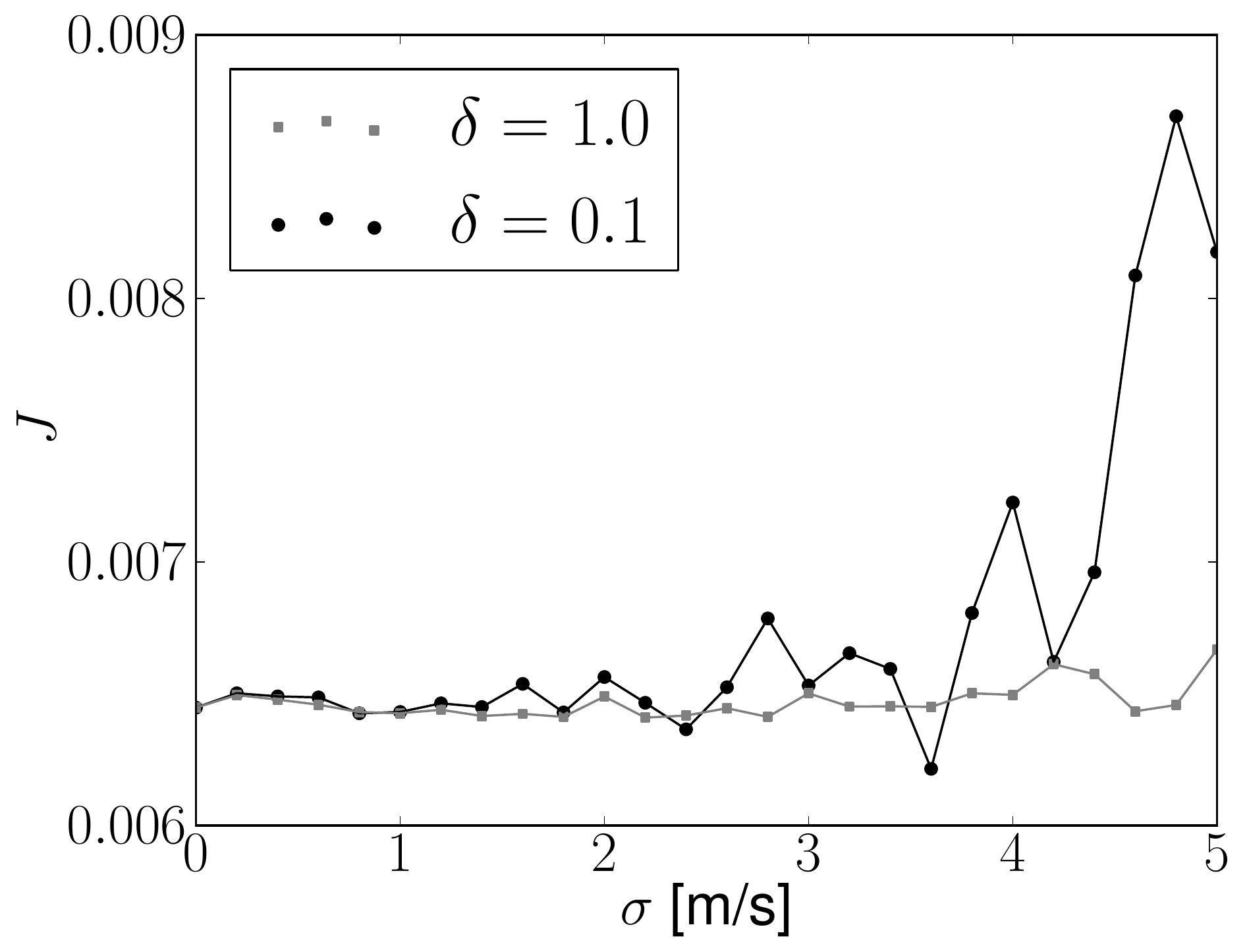}
\caption{Mean linear tracking error as function of turbulence intensity $\sigma$ and correlation rate $\delta$, where $\sigma_{\text{lo}}=\sigma_{\text{la}}=\sigma_{\text{ve}}=\sigma$. Every data point is an
average of 100 simulations.}
\label{fig:performance}

\end{minipage}
\end{figure}

\begin{figure}[h]

\centering

\includegraphics[width=270pt]{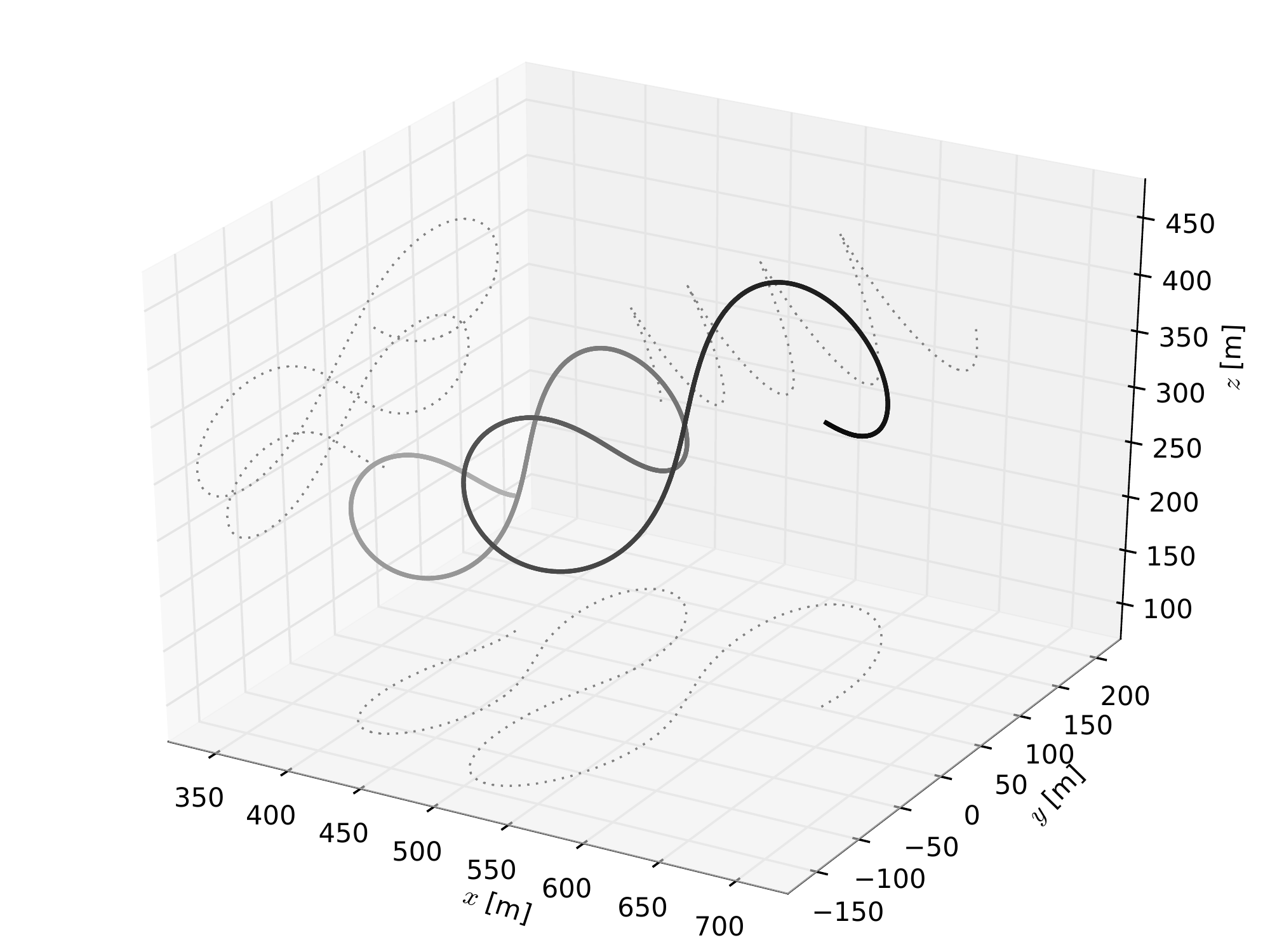}
\caption{Part of a trajectory during tether reel out.}
\label{fig:extending-line}

\end{figure}

Simulation results for a run on constant line length without turbulence are
plotted in Figures \ref{fig:trajectory} and \ref{fig:steering}. Part of a
trajectory with a tether reel out velocity equal to a third of the wind speed
is shown in Figure \ref{fig:extending-line}.

\subsubsection{Turbulence}

In order to quantify the controller performance under the influence of
turbulence we consider the mean linear tracking error
\[
J := \frac{1}{T} \int_{0}^{T} \! \|\boldsymbol\gamma - \boldsymbol\gamma_t\|\, dt.
\]
Figure \ref{fig:performance} plots the mean linear tracking error $J$ for a
single cycle as a function of turbulence
intensity for precise initial control derivative estimates. The turbulence
is modeled according to the simplified Dryden model \cite{Beal93}. The simulation runs
at 100 Hz and every data point is an average of 100 simulations.  The graph
shows that for the investigated range of turbulence intensities the tracking
error remains bounded by a value of $0.009$. Scaling this number with the
distance between the earth tether attachment point and the kite ($100$ m) we obtain a
bound of less than one meter.

\subsubsection{Measurement noise}

\begin{figure}[h]

\centering

\includegraphics[width=160pt]{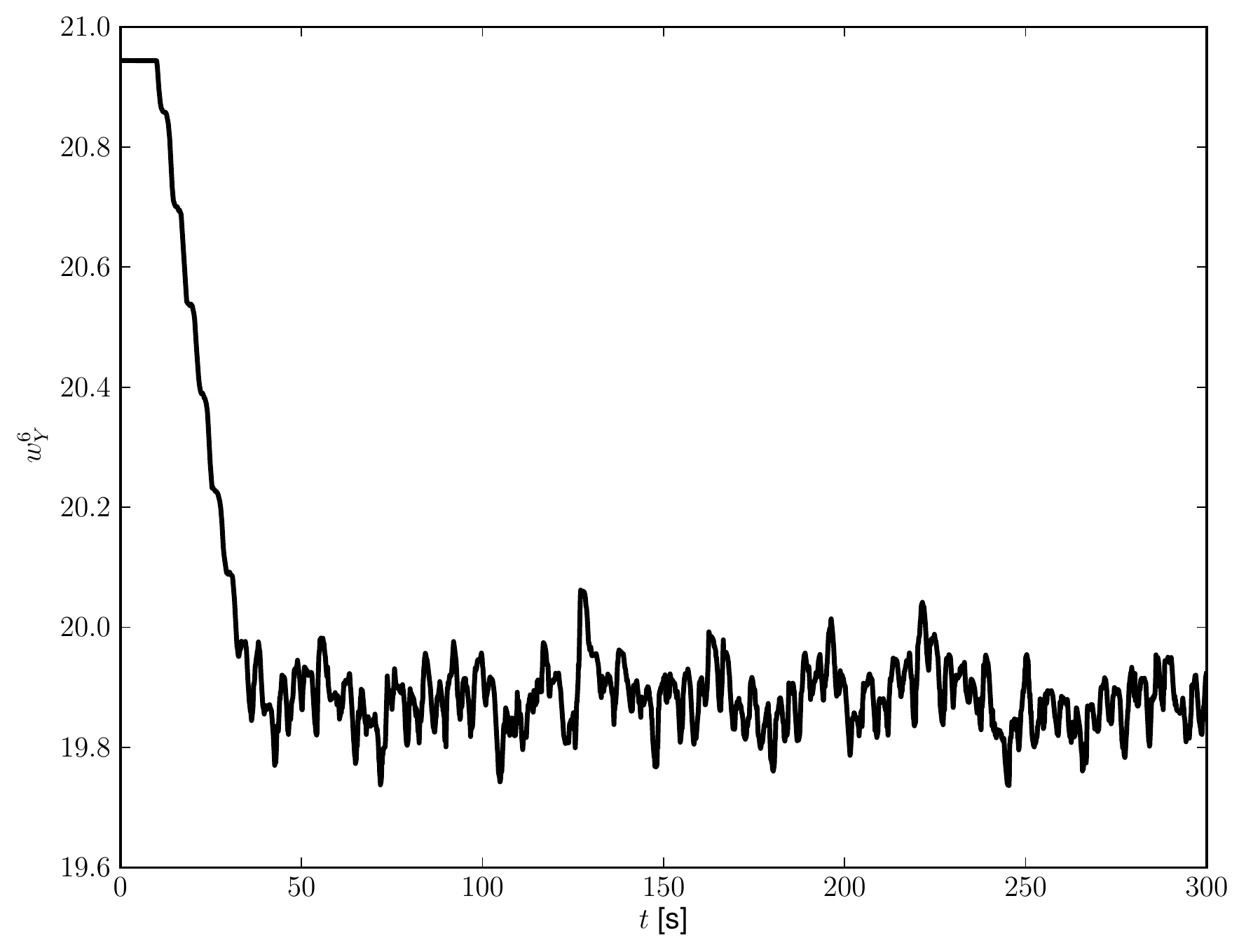}
\caption{Evolution of the $6$th Y-axis control derivative estimate with  Gaussian white noise at the levels given in Table \ref{table:noise-levels}.}
\label{fig:control-deriv-evolution}

\end{figure}

Noisy measurements degrade control performance, which might conceivably lead to problems with adaptivity.
Initial investigations suggest that measurement noise does not induce drift in
the estimates of the control derivatives; see Figure~\ref{fig:control-deriv-evolution} for the evolution of a control derivative estimate
with Gaussian white noise at the levels given in Table~\ref{table:noise-levels}. Assessment of the exact effects of
measurement noise on adaptivity performance is a topic for further investigation and is beyond the scope of this work.

\begin{table}[h]
\centering
\begin{tabular}{|l|l|l|}
\hline
Measurement & Mean & Standard deviation \\
\hline
Position & $0$ m & $2.5$ m \\
Velocity & $0$ m/s & $0.02$ m/s \\
Acceleration & $0$ m/s$^2$ & $0.02$ m/s$^2$ \\
Airspeed & $0$ m/s & $0.9$ m/s \\
Attitude & $0$ deg & $1$ deg \\
\hline
\end{tabular}
\caption{Measurement noise levels.}
\label{table:noise-levels}
\end{table}

\section{Concluding remarks}

In this paper we designed a controller for trajectory tracking control of
kites.  We projected the kite trajectory onto a sphere and used the
differential-geometric notion of turning angle as primary tracking variable.
Based on this dimension reduction we derived an adaptive controller with
minimal modelling requirements.  Repeated simulations with a point-mass model
show our control approach to be robust against turbulence.  A sequel paper will
detail the performance of our controller with multi-body simulations of a
deforming kite on a flexible tether.

\section*{Acknowledgments}

The authors are indebted to Dr. Q. P. Chu for his excellent course on
advanced flight control at Delft University of Technology, to prof. J. H.
van Schuppen, Dr. R. Schmehl and Dr. D. Jeltsema for their comments and
supervision, to Dipl.-Ing. M. Klaus at NTS Energie- und Transportsysteme GmbH for fruitful discussions, and to the entire Laddermill team for their support and
suggestions.

\section*{Appendix A: Geodesic curvature}

In this section we list the main results related to the concept of geodesic curvature referred to
in this paper, adapted from the book by Gray \cite{DifferentialGeometryBook}.
We will not make the notions of surface and curve precise here, but only note
that a \emph{coordinate patch} is a differentiable map between
an open set $\mathcal{U} \subset \mathbb{R}^2$ and $\mathbb{R}^3$: the intuition
is that it associates coordinates to a part of a surface.

In the following, $J$ will denote a counterclockwise rotation by $\pi/2$ in the
tangent plane of the surface under consideration.

We start with the formal definition of the turning angle. Intuitively, the
turning angle is the angle between the tangent of a curve and a reference,
without being restricted to $[-\pi,\pi]$:

\begin{defn}
Let $\mathcal{M}$ be an oriented surface, suppose that $\alpha: (a,b) \to \mathcal{M}$ is a regular curve, and that
$\vec{X}$ is an \emph{everywhere nonzero} vector field along $\boldsymbol\alpha$.
Fix $t_0$ with $a < t_0 < b$. Let $\theta_0$ be a number such that
\[
\frac{\dot{\boldsymbol\alpha}(t_0)}{\|\dot{\boldsymbol\alpha}(t_0)\|} = \cos \theta_0 \frac{\vec{X}(t_0)}{\|\vec{X}(t_0)\|} + \sin \theta_0 \frac{J\vec{X}(t_0)}{\|\vec{X}(t_0)\|}.
\]
Then there exists a unique differentiable function $\theta=\theta[\boldsymbol\alpha,\vec{X}]: (a,b) \to \mathbb{R}$
such that $\theta(t_0)=\theta_0$ and
\[
\frac{\dot{\boldsymbol\alpha}}{\|\dot{\boldsymbol\alpha}\|} = \cos \theta[\boldsymbol\alpha,\vec{X}] \frac{\vec{X}}{\|\vec{X}\|} + \sin \theta[\boldsymbol\alpha,\vec{X}] \frac{J\vec{X}}{\|\vec{X}\|},
\]
at all points on the curve. We call $\theta[\boldsymbol\alpha,\vec{X}]$ the \emph{turning angle} of $\boldsymbol\alpha$
with respect to $\vec{X}$ determined by $\theta_0$ and $t_0$.
\end{defn}

In order to state Liouville's theorem, which allows us to compute the turning
angle of a curve, we need the notion of \emph{geodesic curvature}. Intuitively
speaking the geodesic curvature of a curve measures how far it is from being a
geodesic, that is, the shortest path between two points on the surface:

\begin{defn}
Let $\boldsymbol\alpha: (a,b) \to \mathcal{M}$ be a curve in an oriented regular surface $\mathcal{M}$ in $\mathbb{R}^3$.
Then the \emph{geodesic curvature} of $\boldsymbol\alpha$ is defined
(for $a < t < b$) as
\[
\kappa_g[\boldsymbol\alpha](t)=\frac{\ddot{\boldsymbol\alpha}(t) \cdot J\dot{\boldsymbol\alpha}(t)}{\|\dot{\boldsymbol\alpha}(t)\|^3}.
\]
\end{defn}

Finally, we state Liouville's theorem, which relates the geodesic curvature of a
curve to the geodesic curvatures of the coordinate lines induced by the
coordinate patch and the turning angle of the curve. Let $(\kappa_g)_1$ and
$(\kappa_g)_2$ denote the geodesic curvatures of the coordinate lines $u
\mapsto \vec{x}(u, v)$ and $v \mapsto \vec{x}(u,v)$, where $\vec{x}$ is a given
patch.

\begin{theorem}[Liouville]
Let $\mathcal{M}$ be an oriented surface, and suppose that $\boldsymbol\alpha: (a,b) \to \mathcal{M}$ is a regular curve
whose trace is contained in $\vec{x}(\mathcal{U})$, where $\vec{x}: \mathcal{U} \to \mathcal{M}$ is a coordinate patch
such that $\vec{x}_v=J\vec{x}_u$. Then
\[
\kappa_g[\boldsymbol\alpha] = (\kappa_g)_1 \cos \theta + (\kappa_g)_2 \sin \theta + \frac{\dot{\theta}}{\|\dot{\boldsymbol\alpha}\|},
\]
where $\theta=\theta[\boldsymbol\alpha,\vec{x}_u]$ is the angle between $\boldsymbol\alpha$ and $\vec{x}_u$.
\end{theorem}

\section*{Appendix B: Geodesic distance on the sphere}

In this section we list the main results related to the concept of geodesic distance on the sphere referred to
in this paper, adapted from the report by Bullo, Murray, and Sarti \cite{Bullo95}.

\begin{defn}
The \emph{geodesic distance} between the points $\vec{p}$ and $\vec{q}$ on the unit sphere is defined as
\begin{equation}\label{eq:geodist}
\dist(\vec{p},\vec{q}):=\acos \vec{p} \cdot \vec{q}.
\end{equation}
\end{defn}

\begin{figure}[h]

\centering

\includegraphics[width=120pt]{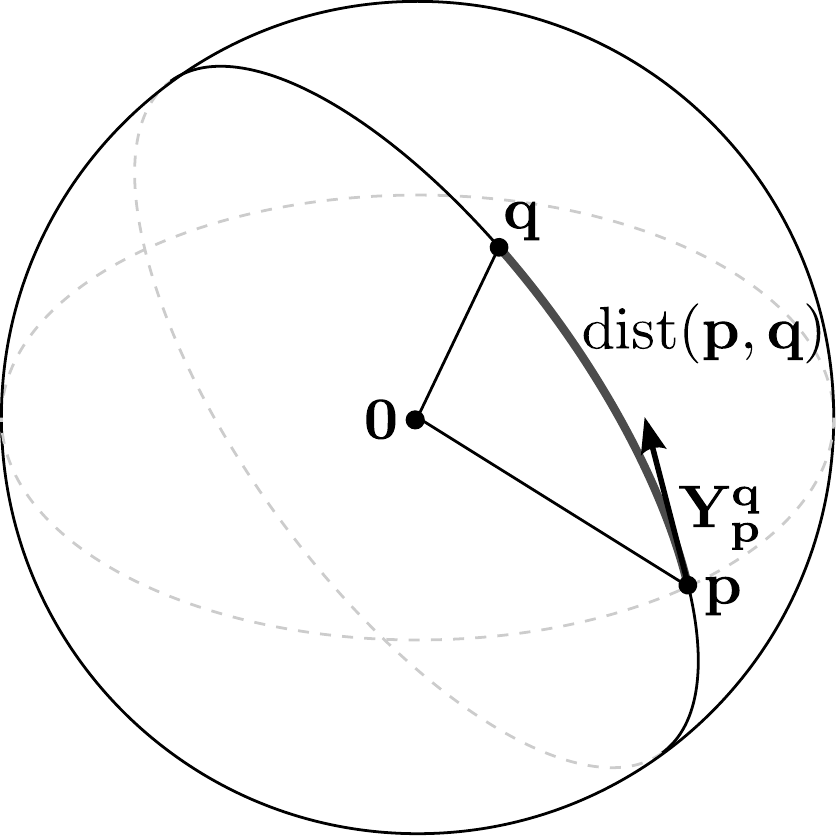}
\caption{Geodesic on the sphere.}
\label{fig:sphere}

\end{figure}

Additionally, provided that $\vec{p} \neq \pm \vec{q}$, there is a uniquely defined unit vector that gives the geodesic
direction at the point $\vec{p}$ towards the point $\vec{q}$:

\begin{defn}
The \emph{geodesic vector} giving the geodesic direction in a point $p$ towards a point $q \neq p$ is defined as
\[
\vec{Y}_{\vec{p}}^{\vec{q}}:= \frac{(\vec{p} \times \vec{q}) \times \vec{p}}{\|(\vec{p} \times \vec{q}) \times \vec{p}\|} = \frac{\vec{q} - (\vec{q} \cdot \vec{p}) \vec{p}}{\|\vec{q} - (\vec{q} \cdot \vec{p}) \vec{p}\|}.
\]
\end{defn}

The time-derivative of the geodesic distance between a trajectory $\vec{p}=\vec{p}(t)$ and a fixed point $\vec{q}$ can be expressed in terms of
the geodesic direction $\vec{Y}_{\vec{p}}^{\vec{q}}$:

\begin{lemma}
Consider a trajectory $\vec{p}=\vec{p}(t)$ on the sphere, such that $\vec{p}(t)$ never passes through the fixed points $\vec{q}$ or $-\vec{q}$.
Then
\[
\frac{d}{dt} \dist(\vec{p}(t),q) = -\dot{\vec{p}} \cdot \vec{Y}_{\vec{p}}^{\vec{q}}.
\]
\end{lemma}
\begin{proof}
Differentiating the definition in Equation (\ref{eq:geodist}), we have
\begin{align*}
\frac{d}{dt} \dist(\vec{p}(t),\vec{q}) & = \frac{d}{dt} \acos \vec{p} \cdot \vec{q} \\
                           & = - \frac{\dot{\vec{p}} \cdot \vec{q}}{\sqrt{1 - (\vec{p} \cdot \vec{q})^2}}
                           = - \frac{\dot{\vec{p}} \cdot (\vec{q} - (\vec{q} \cdot \vec{p})\vec{p})}{\| \vec{p} \times \vec{q}\|}
                           = -\dot{\vec{p}} \cdot \vec{Y}_{\vec{p}}^{\vec{q}}.
\end{align*}
\end{proof}

Let $\operatorname{R}$ be the rotation about $\vec{p} \times \vec{q}$ which maps $\vec{p}$ to $\vec{q}$. We will now extend the previous lemma to cover the time-derivative of the geodesic distance between two trajectories
$\vec{p}=\vec{p}(t)$ and $\vec{q}=\vec{q}(t)$.

\begin{lemma}
\label{lemma:geodist-deriv}
Consider a trajectory $\vec{p}=\vec{p}(t)$ on the sphere, such that $\vec{p}(t)~\neq~\pm \vec{q}(t)$ for all $t > 0$.
Then
\[
\frac{d}{dt} \dist(\vec{p}(t),\vec{q}(t)) = -(\dot{\vec{p}} - \operatorname{R}^T \dot{\vec{q}}) \cdot \vec{Y}_{\vec{p}}^{\vec{q}}.
\]
\end{lemma}
\begin{proof}
By differentiating first with respect to $\vec{p}$ and then with respect to $\vec{q}$, we have
\[
\frac{d}{dt} \dist(\vec{p},\vec{q}) = -\dot{\vec{p}} \cdot \vec{Y}_{\vec{p}}^{\vec{q}} -\dot{\vec{q}} \cdot \vec{Y}_{\vec{q}}^{\vec{p}}.
\]
Since $\operatorname{R}\vec{Y}_{\vec{p}}^{\vec{q}}=-\vec{Y}_{\vec{q}}^{\vec{p}}$, we can simplify the previous expression to
\[
\frac{d}{dt} \dist(\vec{p},\vec{q}) = -(\dot{\vec{p}} - \operatorname{R}^T \dot{\vec{q}}) \cdot \vec{Y}_{\vec{p}}^{\vec{q}}.
\]
\end{proof}

\bibliography{kite-control}
\bibliographystyle{plain}



\end{document}